# Lyapunov exponents for the map that passes through the non-trivial zeros of Riemann zeta-function

J. L. E. da Silva

leonardojade@alu.ufc.br

*Lab. of Quantum Information Technology, Department of Teleinformatic Engineering – Federal University of Ceara - DETI/UFC, C.P. 6007 – Campus do Pici - 60455-970 Fortaleza-Ce, Brazil.*

*Abstract*

The Riemann Hypothesis is the main open problem of Number Theory and several scientists are trying to solve this problem. In this regard, in a recent work [8], a difference equation has been proposed that calculates the nth non-trivial zero in the critical range. In this work, we seek to optimize this estimation by calculating Lyapunov numbers for this non-linear map in order to seek the best value for the bifurcation parameter. Analytical results are presented.

*Key words* – Riemann zeta function. Lambert *W* function. Lyapunov exponents.

## 1.Introduction

The Riemann $\zeta(z)$ function is a central point in the study of the Riemann Hypothesis [1]. This function has a close relationship with the prime numbers *P*

$$\zeta(z) = \sum_{n=1}^{\infty} n^{-z} = \prod_{p \in P}(1 - p^{-z})^{-1} . \qquad (1)$$

Its analytical continuity can be obtained through

$$\zeta(z) = \left(1 - 2^{1-z}\right)^{-1} \sum_{n=1}^{\infty} (-1)^{n-1} n^{-z} , \qquad (2)$$

for $0 \leq Re(z) < 1$ in the critical range.

## 2. Statistics about the non-trivial zeros of the function $\zeta(z)$

Statistics in the critical range of zeros of the Riemann function is one of the possibilities for current research on the Riemann Hypothesis [2]. All non-trivial zeros of $\zeta(z)$ are known to be complex numbers with positive imaginary part in the form $\rho = 1/2 + it_n$ with $0 < t_n \leq t$. The number of non-trivial zeros with multiplicity of the function up to a range t can be calculated using [3]

$$N(t) = \frac{t}{2\pi}\log\left(\frac{t}{2\pi}\right) - \frac{t}{2\pi} + \frac{7}{8} + S(t) + O(t^{-1}), \tag{3}$$

where $S(t)$ represents the argument function [4]

$$S(t) = \frac{1}{\pi}\arg\zeta\left(\frac{1}{2} + it\right). \tag{4}$$

In [5], França and Leclair investigated the non-trivial zeros statistic through the solution of transcendental equation

$$\frac{t_n}{2\pi}\log\left(\frac{t_n}{2\pi e}\right) + \lim_{\delta \to 0^+}\frac{1}{\pi}\arg\zeta\left(\frac{1}{2} + \delta + it_n\right) = n - \frac{11}{8}, \tag{5}$$

for $(n = 1, 2, \ldots, t)$. They obtained the following estimate for $t_n$

$$\tilde{t}_n = \frac{2\pi\left(n - \frac{11}{8}\right)}{W_0\left[e^{-1}\left(n - \frac{11}{8}\right)\right]}, \tag{6}$$

where $W_0$ is the main branch of the Lambert $W$ function [6,7]. Based on this proposal, Ramos and Silva [8] proposed a difference equation that calculates the $n$-th non-trivial zero of the function $\zeta(z)$ in the critical line

$$t_n^k = \frac{2\pi\left(n - \frac{11}{8} - \frac{\delta}{\pi}\arg\zeta\left(\frac{1}{2} + it_n^{k-1}\right)\right)}{W_0\left(\frac{1}{e}\left(n - \frac{11}{8} - \frac{\delta}{\pi}\arg\zeta\left(\frac{1}{2} + it_n^{k-1}\right)\right)\right)}. \tag{7}$$

## 3. Lyapunov numbers for the nonlinear map that calculates the non-trivial zeros of the function $\zeta(z)$ on the critical line

One way to search for the best value of $\delta$ de (7), considering that $\delta = 1$ is not always the best value for the bifurcation parameter of the map, in this perspective, we seek to calculate the Lyapunov numbers for the map of (7) [9]. Considering $t_n^k = t$ and deriving the map as a function of $t$ we obtain as a result

$$\Delta(\delta) = \lim_{N\to\infty} \frac{1}{N} \sum_{n=0}^{N} \ln\left|\frac{\left[\left(\delta\log\frac{t}{2\pi} - 2\pi\delta O(t^{-2})\right)\right]\left[W_0\left(\frac{1}{e}\left(n - \frac{11}{8} - \frac{\delta}{\pi}\arg\zeta\left(\frac{1}{2} + it\right)\right)\right) - \frac{e^{-W_0\left(\frac{1}{e}\left(n - \frac{11}{8} - \frac{\delta}{\pi}\arg\zeta\left(\frac{1}{2} + it\right)\right)\right)}}{1 + W_0\left(\frac{1}{e}\left(n - \frac{11}{8} - \frac{\delta}{\pi}\arg\zeta\left(\frac{1}{2} + it\right)\right)\right)}\right]\left(\frac{1}{e}\left(n - \frac{11}{8}\frac{\delta}{\pi}\arg\zeta\left(\frac{1}{2} + it\right)\right)\right)}{\left(W_0\left(\frac{1}{e}\left(n - \frac{11}{8} - \frac{\delta}{\pi}\arg\zeta\left(\frac{1}{2} + it\right)\right)\right)\right)^2}\right|$$

(8)

whose numerical calculation must consider that $O(t^{-2})$ is equal to zero.

## 4. Conclusions

Using (7) to numerically estimate the numbers of Lyapunov (8), we can find the bifurcation points that indicate the best $\delta$ values to be used to obtain with better precision the complex part of the non-trivial zeros of the function $\zeta(z)$.

**Appendix**

$$\frac{d}{dt}\left[\frac{2\pi\left(n-\frac{11}{8}-\frac{\delta}{\pi}\arg\zeta\left(\frac{1}{2}+it\right)\right)}{W_0\left(\frac{1}{e}\left(n-\frac{11}{8}-\frac{\delta}{\pi}\arg\zeta\left(\frac{1}{2}+it\right)\right)\right)}\right]=$$

Let $S(t)=\frac{1}{\pi}\arg\zeta\left(\frac{1}{2}+it\right)$, we have to

$$\frac{d}{dt}\left[\frac{2\pi\left(n-\frac{11}{8}-\delta S(t)\right)}{W_0\left(\frac{1}{e}\left(n-\frac{11}{8}-\delta S(t)\right)\right)}\right]=$$

$$\left[\frac{\frac{d}{dt}\left(2\pi\left(n-\frac{11}{8}-\delta S(t)\right)\right)\left(W_0\left(\frac{1}{e}\left(n-\frac{11}{8}-\delta S(t)\right)\right)\right)-\frac{d}{dt}\left(W_0\left(\frac{1}{e}\left(n-\frac{11}{8}-\delta S(t)\right)\right)\right)\left(2\pi\left(n-\frac{11}{8}-\delta S(t)\right)\right)}{\left(W_0\left(\frac{1}{e}\left(n-\frac{11}{8}-\delta S(t)\right)\right)\right)^2}\right]$$

$=$

by [4,6], we have to

$$\frac{d}{dt}(S(t))=-\frac{1}{2\pi}\log\frac{t}{2\pi}+O(t^{-2}),$$

e

$$\frac{d}{dt}(W_0(u(t)))=\left(\frac{e^{-W_0(u(t))}}{1+W_0(u(t))}\right)\frac{d}{dt}(u(t)),$$

finding

$$\frac{\left(-2\pi\delta\frac{d}{dt}(S(t))\right)\left(W_0\left(\frac{1}{e}\left(n-\frac{11}{8}-\delta S(t)\right)\right)\right)-\left(\frac{e^{-W_0\left(\frac{1}{e}\left(n-\frac{11}{8}-\delta S(t)\right)\right)}}{1+W_0\left(\frac{1}{e}\left(n-\frac{11}{8}-\delta S(t)\right)\right)}\right)\frac{d}{dt}\left(\frac{1}{e}\left(n-\frac{11}{8}-\delta S(t)\right)\right)\left(2\pi\left(n-\frac{11}{8}-\delta S(t)\right)\right)}{\left(W_0\left(\frac{1}{e}\left(n-\frac{11}{8}-\delta S(t)\right)\right)\right)^2}=$$

$$\frac{\left(-2\pi\delta\frac{d}{dt}(S(t))\right)\left(W_0\left(\frac{1}{e}\left(n-\frac{11}{8}-\delta S(t)\right)\right)\right)-\left(\frac{e^{-W_0\left(\frac{1}{e}\left(n-\frac{11}{8}-\delta S(t)\right)\right)}}{1+W_0\left(\frac{1}{e}\left(n-\frac{11}{8}-\delta S(t)\right)\right)}\right)\left(-\frac{\delta}{e}\frac{d}{dt}S(t)\right)\left(2\pi\left(n-\frac{11}{8}-\delta S(t)\right)\right)}{\left(W_0\left(\frac{1}{e}\left(n-\frac{11}{8}-\delta S(t)\right)\right)\right)^2}=$$

$$\frac{\left(-2\pi\delta\frac{d}{dt}(S(t))\right)\left(W_0\left(\frac{1}{e}\left(n-\frac{11}{8}-\delta S(t)\right)\right)\right)-\left(\frac{e^{-W_0\left(\frac{1}{e}\left(n-\frac{11}{8}-\delta S(t)\right)\right)}}{1+W_0\left(\frac{1}{e}\left(n-\frac{11}{8}-\delta S(t)\right)\right)}\right)\left(-\frac{2\pi\delta}{e}\frac{d}{dt}S(t)\right)\left(\left(n-\frac{11}{8}-\delta S(t)\right)\right)}{\left(W_0\left(\frac{1}{e}\left(n-\frac{11}{8}-\delta S(t)\right)\right)\right)^2}=$$

$$\frac{\left(-2\pi\delta\frac{d}{dt}(S(t))\right)\left[\left(W_0\left(\frac{1}{e}\left(n-\frac{11}{8}-\delta S(t)\right)\right)\right)-\left(\frac{e^{-W_0\left(\frac{1}{e}\left(n-\frac{11}{8}-\delta S(t)\right)\right)}}{1+W_0\left(\frac{1}{e}\left(n-\frac{11}{8}-\delta S(t)\right)\right)}\right)\left(\frac{1}{e}\left(n-\frac{11}{8}-\delta S(t)\right)\right)\right]}{\left(W_0\left(\frac{1}{e}\left(n-\frac{11}{8}-\delta S(t)\right)\right)\right)^2}=$$

$$\frac{\left(-2\pi\delta\left(-\frac{1}{2\pi}\log\frac{t}{2\pi}+O(t^{-2})\right)\right)\left[\left(W_0\left(\frac{1}{e}\left(n-\frac{11}{8}-\delta S(t)\right)\right)\right)-\left(\frac{e^{-W_0\left(\frac{1}{e}\left(n-\frac{11}{8}-\delta S(t)\right)\right)}}{1+W_0\left(\frac{1}{e}\left(n-\frac{11}{8}-\delta S(t)\right)\right)}\right)\left(\frac{1}{e}\left(n-\frac{11}{8}-\delta S(t)\right)\right)\right]}{\left(W_0\left(\frac{1}{e}\left(n-\frac{11}{8}-\delta S(t)\right)\right)\right)^2}=$$

$$\frac{\left(\left(\delta\log\frac{t}{2\pi}-2\pi\delta O(t^{-2})\right)\right)\left[\left(W_0\left(\frac{1}{e}\left(n-\frac{11}{8}-\delta S(t)\right)\right)\right)-\left(\frac{e^{-W_0\left(\frac{1}{e}\left(n-\frac{11}{8}-\delta S(t)\right)\right)}}{1+W_0\left(\frac{1}{e}\left(n-\frac{11}{8}-\delta S(t)\right)\right)}\right)\left(\frac{1}{e}\left(n-\frac{11}{8}-\delta S(t)\right)\right)\right]}{\left(W_0\left(\frac{1}{e}\left(n-\frac{11}{8}-\delta S(t)\right)\right)\right)^2}=$$

$$\Delta(\delta) = \lim_{N\to\infty} \frac{1}{N} \sum_{n=0}^{N} \ln \left| \frac{\left[\left(\delta\log\frac{t}{2\pi} - 2\pi\delta O(t^{-2})\right)\right]\left[\left(W_0\left(\frac{1}{e}\left(n - \frac{11}{8} - \delta S(t)\right)\right)\right) - \left(\frac{e^{-W_0\left(\frac{1}{e}\left(n - \frac{11}{8} - \delta S(t)\right)\right)}}{1 + W_0\left(\frac{1}{e}\left(n - \frac{11}{8} - \delta S(t)\right)\right)}\right)\left(\frac{1}{e}\left(n - \frac{11}{8} - \delta S(t)\right)\right)\right]}{\left(W_0\left(\frac{1}{e}\left(n - \frac{11}{8} - \delta S(t)\right)\right)\right)^2} \right|$$

=

$$\Delta(\delta) = \lim_{N\to\infty} \frac{1}{N} \sum_{n=0}^{N} \ln \left| \frac{\left[\left(\delta\log\frac{t}{2\pi} - 2\pi\delta O(t^{-2})\right)\right]\left[\left(W_0\left(\frac{1}{e}\left(n - \frac{11}{8} - \frac{\delta}{\pi}\arg\zeta\left(\frac{1}{2} + it\right)\right)\right)\right) - \left(\frac{e^{-W_0\left(\frac{1}{e}\left(n - \frac{11}{8} - \frac{\delta}{\pi}\arg\zeta\left(\frac{1}{2} + it\right)\right)\right)}}{1 + W_0\left(\frac{1}{e}\left(n - \frac{11}{8} - \frac{\delta}{\pi}\arg\zeta\left(\frac{1}{2} + it\right)\right)\right)}\right)\left(\frac{1}{e}\left(n - \frac{11}{8}\frac{\delta}{\pi}\arg\zeta\left(\frac{1}{2} + it\right)\right)\right)\right]}{\left(W_0\left(\frac{1}{e}\left(n - \frac{11}{8} - \frac{\delta}{\pi}\arg\zeta\left(\frac{1}{2} + it\right)\right)\right)\right)^2} \right|$$

.

## References

[1] CONREY, J. The riemann hypothesis. *Notices of the AMS*, 2003, 50.3: 341-353.

[2] CHOWLA, S. *The Riemann hypothesis and Hilbert's tenth problem*. CRC Press, 1987.

[3] CARNEIRO, E. , et al. Bandlimited approximations and estimates for the Riemann zeta-function. *Publicacions Matemàtiques*, 2019, 63.2: 601-661.

[4] KARATSUBA, A., et al. The argument of the Riemann zeta function. *Russian Mathematical Surveys*, 2005, 60.3: 433-488.

[5] FRANÇA, G.; LECLAIR, A. Statistical and other properties of Riemann zeros based on an explicit equation for the *n*-th zero on the critical line. *arXiv preprint arXiv:1307.8395*, 2013.

[6] R. M. Corless, G. H. Gonnet, D. E. G. Hare, D. J. Jeffrey and D. E. Knuth, On the Lambert W function, Advances in Computational Mathematics, vol. 5, 329 – 359, 1996.

[7] Valluri, S. R., Jeffrey, D. J., & Corless, R. M. Some applications of the Lambert W function to physics. *Canadian Journal of Physics*, 78(9), 823-831. (2000).

[8] DA SILVA, G. B.; RAMOS, R. V. A Non-Linear Difference Equation for Calculation of the Zeros of the Riemann Zeta-Function on the Critical Line. *arXiv preprint arXiv:1810.01823*, 2018.

[9] DE SOUZA, S. CALDAS, I. Calculation of Lyapunov exponents in systems with impacts. *Chaos, Solitons & Fractals*, 19.3: 569-579. 2004.